\documentclass[10pt]{article}
\usepackage{lmodern}
\usepackage{amsmath}
\usepackage[T1]{fontenc}
\usepackage[utf8]{inputenc}
\usepackage{authblk}
\usepackage{amsfonts}
\usepackage{graphicx}
\usepackage{tkz-euclide}
\usetkzobj{all}
\tikzset{elegant/.style={smooth,thick,samples=50,cyan}}
\tikzset{eaxis/.style={->,>=stealth}}
\usepackage{color,tikz}
\usetikzlibrary{decorations.pathreplacing,calc}
\usepackage{rotating}
\usepackage{amssymb}
\usepackage[english]{babel}
\usepackage{color}
\usepackage{amsthm}
\usepackage{graphicx}
\usepackage{mathrsfs}
\usepackage{makecell}
\usepackage{microtype}
\usepackage{mathscinet}
\usepackage{array}
\usepackage{multirow}
\usepackage{enumerate}
\usepackage[cal=boondoxo,bb=ams]{mathalfa}
\usepackage{hyperref}
\hypersetup{hidelinks}

\newtheorem{defn}{Definition}[section]
\newtheorem{theorem}{Theorem}[section]
\newtheorem{prop}{Proposition}[section]

\newtheorem{remark}{Remark}[section]

\newcommand{\ml}{\mathcal}
\newcommand{\mb}{\mathbb}

\title{Blow-up of solutions to Nakao's problem via an iteration argument}
\author[1]{Wenhui Chen\thanks{Corresponding author: Wenhui Chen (wenhui.chen.math@gmail.com)}}
\author[1]{Michael Reissig}
\affil[1]{Institute of Applied Analysis, Faculty of Mathematics and Computer Science, Technical University Bergakademie Freiberg, Pr\"uferstra{\ss}e 9 09596 Freiberg, Germany}
\date{}
\begin{document}

		\maketitle
\begin{abstract}
In this paper, we consider blow-up behavior of weak solutions to a weakly coupled system for a semilinear damped wave equation and a semilinear wave equation in $\mb{R}^n$. This problem is part of the  so-called Nakao's problem proposed by Professor Mitsuhiro Nakao (Kyushu university) for a critical relation between the exponents $p$ and $q$. By applying an iteration method for unbounded multipliers with a slicing procedure, we prove blow-up of weak solutions for Nakao's problem even for small data. We improve the blow-up result and upper bound estimates for lifespan comparing with the previous research, especially, in higher dimensional cases.
	\medskip\\
	\textbf{AMS Classification (2010)} Primary 35L52; Secondary 35B44
\end{abstract}
\fontsize{12}{15}
\selectfont

\section{Introduction} \label{Sec1}
We investigate blow-up of solutions in finite time of weak solutions to the Cauchy problem for a weakly coupled system of a semilinear damped wave equation and a semilinear wave equation, namely, for
weak solutions to
\begin{equation}\label{Nakao's problem}
\begin{cases}
u_{tt}-\Delta u+u_t=|v|^p,&x\in\mb{R}^n, \ t>0,\\
v_{tt}-\Delta v=|u|^q,&x\in\mb{R}^n,\ t>0,\\
(u,u_t,v,v_t)(0,x)=\varepsilon(u_0,u_1,v_0,v_1)(x),&x\in\mb{R}^n,
\end{cases}
\end{equation}
where $p,q>1$ and $\varepsilon$ is a positive parameter describing the size of initial data. The problem of critical curve of exponents $p$ and $q$ for the weakly coupled system \eqref{Nakao's problem} was proposed by Professor Mitsuhiro Nakao, Emeritus of Kyushu University (see also \cite{NishiharaWakasugi2015,Wakasugi2017}). Here, ``critical curve" stands for the threshold condition of a pair of exponents $(p,q)$ between global (in time) existence of small data weak solutions and blow-up of weak solutions even for small data. Recently, by employing the test function method (see, for example, \cite{MP01book,Zhang01}) the author of \cite{Wakasugi2017} proved that if the following conditions:
\begin{align}\label{Condition Wakasugi}
\alpha_{\mathrm{N,W}}:=\max\left\{\frac{q/2+1}{pq-1}+\frac{1}{2},\frac{q+1}{pq-1},\frac{p+1}{pq-1}\right\}\geqslant \frac{n}{2},
\end{align}
and
\begin{align*}
1<p,q<\infty \ (n=1,2),\ \ 1<p,q\leqslant \frac{n}{n-2}\ (n\geqslant 3),
\end{align*}
are satisfied, then every local (in time) weak solution $(u,v)$ to the weakly coupled system \eqref{Nakao's problem} blows up in finite time. Nevertheless, the curve $\alpha_{\text{N,W}}=n/2$ in the $p-q$ plane for a pair of exponents $(p,q)$ seems optimal only when $n=1$. Precisely, for the case $n=1$, under the condition $\alpha_{\mathrm{N,W}}\geqslant1/2$ (or, equivalently, $1<p,q<\infty$) every local (in time) solution blows up. But, in the case $n\geqslant 2$, the condition \eqref{Condition Wakasugi} seems not to be optimal. The main goal in this paper is to improve the blow-up result stated in \cite{Wakasugi2017}, particularly, in higher dimensional cases.

We sketch now some historical background related to the weakly coupled system \eqref{Nakao's problem}. Since Nakao's problem \eqref{Nakao's problem} is in some sense related to weakly coupled systems of semilinear damped wave equations and semilinear wave equations, we would like to recall some results for these systems, respectively.

On the one hand, the weakly coupled system of semilinear wave equations
\begin{equation}\label{weaklycoupledwave}
\begin{cases}
u_{tt}-\Delta u=|v|^{p},&x\in\mathbb{R}^n,\ t>0,\\
v_{tt}-\Delta v=|u|^{q},&x\in\mathbb{R}^n,\ t>0,\\
(u,u_t,v,v_t)(0,x)=(u_0,u_1,v_0,v_1)(x),&x\in\mathbb{R}^n,
\end{cases}
\end{equation}
for $n\geqslant1$ with $p,q>1$, has been widely studied in recent years. The papers \cite{DelS97,DGM,DM,AKT00,KT03,Kur05,GTZ06,KTW12} investigated that the critical curve for the weakly coupled system \eqref{weaklycoupledwave} is described by the condition
\begin{equation}\label{Critical curve wave}
\alpha_{\mathrm{W}}:=\max\left\{\frac{p+2+q^{-1}}{pq-1},\frac{q+2+p^{-1}}{pq-1}\right\}=\frac{n-1}{2}.
\end{equation}
In other words, if $\alpha_{\mathrm{W}}<(n-1)/2$, then there exists a unique global (in time) weak solution for small data. On the contrary, if $\alpha_{\mathrm{W}}\geqslant (n-1)/2$, in general, local (in time) weak solutions blow up. Especially, we should underline that the approach for proving blow-up results for \eqref{weaklycoupledwave} is mainly based on a generalized Kato's type lemma or an iteration argument.

On the other hand, let us turn to recall some results for the weakly coupled system of semilinear classical damped wave equations
\begin{equation}\label{weaklycoupleddamped}
\begin{cases}
u_{tt}-\Delta u+u_t=|v|^{p},&x\in\mathbb{R}^n,\ t>0,\\
v_{tt}-\Delta v+v_t=|u|^{q},&x\in\mathbb{R}^n,\ t>0,\\
(u,u_t,v,v_t)(0,x)=(u_0,u_1,v_0,v_1)(x),&x\in\mathbb{R}^n,
\end{cases}
\end{equation}
for $n\geqslant 1$ with $p,q>1$. The critical curve for the weakly coupled system \eqref{weaklycoupleddamped} is characterized by the condition
\begin{equation}\label{Critical curve damped wave}
\alpha_{\mathrm{DW}}:=
\max\left\{\frac{p+1}{pq-1},\frac{q+1}{pq-1}\right\}=\frac{n}{2},
\end{equation}
which has been investigated by the authors of \cite{SunWang2007,Narazaki2009,Nishi12,NishiharaWakasugi}. To derive nonexistence results for global (in time) weak solutions to the weakly coupled system \eqref{weaklycoupleddamped}, the authors applied the test function method, which is useful to deal with semilinear classical damped wave models with effective damping.

From the above results of critical curves for the weakly coupled systems \eqref{weaklycoupledwave} and \eqref{weaklycoupleddamped}, we may expect that the critical curve for Nakao's problem \eqref{Nakao's problem} is influenced by the relations \eqref{Critical curve wave} and \eqref{Critical curve damped wave}. However, we should underline that the critical curve for the weakly coupled system \eqref{Nakao's problem} is not a simple combination of \eqref{Critical curve wave} and \eqref{Critical curve damped wave} because the critical curve to Nakao's problem \eqref{Nakao's problem} seems to be influenced by varying degrees between semilinear wave equations and semilinear damped wave equations. Let us now focus on the question for blow-up of weak solutions. From the technical point of view, quite different methods are applied to derive blow-up results for the weakly coupled systems \eqref{weaklycoupledwave} and \eqref{weaklycoupleddamped}. Therefore, it is significant for us to find a suitably unified approach to deal with Nakao's problem \eqref{Nakao's problem}. As mentioned above, \cite{Wakasugi2017} proved blow-up of weak solutions for Nakao's problem \eqref{Nakao's problem} by using the test function method. However, for higher dimensional cases ($n\geqslant 4$) the condition fulfills \begin{align*}
\alpha_{\mathrm{N,W}}\geqslant\frac{n-1}{2}\,\,\, \mbox{iff}\,\,\, \  \alpha_{\mathrm{DW}}\geqslant \frac{n-1}{2}.
\end{align*}
In other words, for higher dimensional cases, the blow-up result from \cite{Wakasugi2017} is the same as those for weakly coupled system of classical damped wave equations.
In this paper, we mainly improve the blow-up result from \cite{Wakasugi2017} in higher dimensional cases $n\geqslant4$ by using an iteration method. With the aim of adapting the semilinear damped wave equation, we employ such method for unbounded multipliers, which is caused by the friction $u_t$. More precisely, motivated by \cite{AKT00,ChenPalmieri201901}, we propose a slicing procedure of the domain of integration by suitable sequences (see \eqref{Seq_Nakao} and \eqref{Seq_Nakao_2} later). Finally, thanks to the iteration argument, the blow-up result from \cite{Wakasugi2017} is partially improved for $2\leqslant n\leqslant 3$ and completely improved for $n\geqslant 4$. Simultaneously, some estimates for upper bounds of lifespan can be derived.

\section{Main results} \label{Sec2}
\setcounter{equation}{0}
Before stating the main results of this paper, let us introduce a suitable notion of energy solutions to the weakly coupled system \eqref{Nakao's problem}.
\begin{defn}\label{DefnEnergySolution}
	Let $(u_0,u_1,v_0,v_1)\in \left(H^1(\mb{R}^n)\times L^2(\mb{R}^n)\right)\times \left(H^1(\mb{R}^n)\times L^2(\mb{R}^n)\right)$. We say that $(u,v)$ is an energy solution of the weakly coupled system \eqref{Nakao's problem} on $[0,T)$ if
	\begin{align*}
	u\in\ml{C}\left([0,T),H^1(\mb{R}^n)\right)\cap \ml{C}^1\left([0,T),L^2(\mb{R}^n)\right)\cap L_{\mathrm{loc}}^q([0,T)\times\mb{R}^n),\\
	v\in\ml{C}\left([0,T),H^1(\mb{R}^n)\right)\cap \ml{C}^1\left([0,T),L^2(\mb{R}^n)\right)\cap L_{\mathrm{loc}}^p([0,T)\times\mb{R}^n),
	\end{align*}
	satisfies $(u,v)(0,\cdot)=(u_0,v_0)$ in $H^1(\mb{R}^n)\times H^1(\mb{R}^n)$ and the following integral relations:
	\begin{align}\label{Defn.Energy.Solution 01}
	&\int_0^t\int_{\mb{R}^n}\left(-u_t(s,x)\,\phi_s(s,x)+u_t(s,x)\,\phi(s,x)+\nabla u(s,x)\cdot\nabla\phi(s,x)\right)\mathrm{d}x\,\mathrm{d}s\notag\\
	&+\int_{\mb{R}^n}u_t(t,x)\,\phi(t,x)\,\mathrm{d}x-\int_{\mb{R}^n}u_1(x)\,\phi(0,x)\,\mathrm{d}x=\int_0^t\int_{\mb{R}^n}|v(s,x)|^p\,\phi(s,x)\,\mathrm{d}x\,\mathrm{d}s
	\end{align}
	and
	\begin{align}\label{Defn.Energy.Solution 02}
	&\int_0^t\int_{\mb{R}^n}\left(-v_t(s,x)\,\psi_s(s,x)+\nabla v(s,x)\cdot\nabla\psi(s,x)\right)\mathrm{d}x\,\mathrm{d}s\notag\\
	&+\int_{\mb{R}^n}v_t(t,x)\,\psi(t,x)\,\mathrm{d}x-\int_{\mb{R}^n}v_1(x)\,\psi(0,x)\,\mathrm{d}x=\int_0^t\int_{\mb{R}^n}|u(s,x)|^q\,\psi(s,x)\,\mathrm{d}x\,\mathrm{d}s
	\end{align}
	for any test functions $\phi,\psi\in\ml{C}_0^{\infty}([0,T)\times\mb{R}^n)$ and any $t\in(0,T)$.
\end{defn}
Clearly, the application of further steps of integration by parts in \eqref{Defn.Energy.Solution 01} and \eqref{Defn.Energy.Solution 02}, we obtain
\begin{align*}
&\int_0^t\int_{\mb{R}^n}u(s,x)\left(\phi_{ss}(s,x)-\Delta\phi(s,x)-\phi_s(s,x)\right)\mathrm{d}x\,\mathrm{d}s\\
&+\int_{\mb{R}^n}\left(u_t(t,x)\,\phi(t,x)+u(t,x)\,\phi(t,x)-u(t,x)\,\phi_s(t,x)\right)\mathrm{d}x\\
&-\int_{\mb{R}^n}\left(u_1(x)\,\phi(0,x)+u_0(x)\,\phi(0,x)-u_0(x)\,\phi_s(0,x)\right)\mathrm{d}x=\int_0^t\int_{\mb{R}^n}|v(s,x)|^p\,\phi(s,x)\,\mathrm{d}x\,\mathrm{d}s
\end{align*}
and
\begin{align*}
&\int_0^t\int_{\mb{R}^n}v(s,x)\left(\psi_{ss}(s,x)-\Delta\psi(s,x)\right)\mathrm{d}x\,\mathrm{d}s+\int_{\mb{R}^n}\left(v_t(t,x)\,\psi(t,x)-v(t,x)\,\psi_s(t,x)\right)\mathrm{d}x\\
&-\int_{\mb{R}^n}\left(v_1(x)\,\psi(0,x)-v_0(x)\,\psi_s(0,x)\right)\mathrm{d}x=\int_0^t\int_{\mb{R}^n}|u(s,x)|^q\,\psi(s,x)\,\mathrm{d}x\,\mathrm{d}s.
\end{align*}
Letting $t\rightarrow T$, we find that $(u,v)$ fulfills the definition of weak solutions to Nakao's problem \eqref{Nakao's problem}.
\medskip

Let us state two propositions for blow-up of energy solutions to Nakao's problem \eqref{Nakao's problem}. They are based on different lower bound estimates for functionals appearing in the blow-up dynamic.
\begin{prop}\label{Thm.Main.Blow-Up}
	Let us consider $p,q>1$ if $n=1,2$, and $1<p,q\leqslant n/(n-2)$ if $n\geqslant 3$ such that
	\begin{align}\label{Assumption.Condition}
	\alpha_{0}:=\max\left\{\frac{q/2+1}{pq-1},\frac{2+p^{-1}}{pq-1}\right\}>\frac{n-1}{2}.
	\end{align}
	Let us assume that $(u_0,u_1,v_0,v_1)\in \left(H^1(\mb{R}^n)\times L^2(\mb{R}^n)\right)\times \left(H^1(\mb{R}^n)\times L^2(\mb{R}^n)\right)$ are nonnegative and compactly supported functions with supports contained in $B_R$ for some $R>0$ such that $u_0,v_1$ are not identically zero. Let $(u,v)$ be the local (in time) energy solution to Nakao's problem \eqref{Nakao's problem} according to Definition \ref{DefnEnergySolution} with lifespans $T=T(\varepsilon)$. Then these solutions satisfy
	\begin{align}\label{Supp_uv}
	\mathrm{supp}\, u,\,\mathrm{supp}\,v\subset \{(t,x)\in[0,T)\times\mb{R}^n:|x|\leqslant R+t\}.
	\end{align}
	Moreover, there exists a positive constant $\varepsilon_0=\varepsilon_0(u_0,u_1,v_0,v_1,n,p,q,R)$ such that for any $\varepsilon\in(0,\varepsilon_0]$ the solution $(u,v)$ blows up in finite time. Furthermore, the upper bound estimate for the lifespans
	\begin{align*}
	T(\varepsilon)\leqslant C\varepsilon^{-1/\max\{F_1(n,p,q),F_2(n,p,q)\}}
	\end{align*}
	holds, where $C>0$ is a constant independent of $\varepsilon$ and
	\begin{align}
	F_1(n,p,q)&:= \frac{2+p^{-1}}{pq-1}-\frac{n-1}{2},\label{Defn.F1}\\
	F_2(n,p,q)&:= \frac{1+2q^{-1}}{pq-1}-\frac{n-1}{q}.\label{Defn.F2}
	\end{align}
\end{prop}
\begin{prop}\label{Thm.Main.Blow-Up_2}
	Let us consider $p,q>1$ if $n=1,2$, and $1<p,q\leqslant n/(n-2)$ if $n\geqslant 3$ such that
	\begin{align}\label{Assumption.Condition_2}
	\alpha_{1}:=\max\left\{\frac{q/2+1}{pq-1},\frac{1/2+p}{pq-1}-\frac{1}{2}\right\}>\frac{n-1}{2}.
	\end{align}
	Let us assume that $(u_0,u_1,v_0,v_1)\in \left(H^1(\mb{R}^n)\times L^2(\mb{R}^n)\right)\times \left(H^1(\mb{R}^n)\times L^2(\mb{R}^n)\right)$ have the same assumption as in Theorem \ref{Thm.Main.Blow-Up}. Let $(u,v)$ be the energy solution to Nakao's problem \eqref{Nakao's problem} according to Definition \ref{DefnEnergySolution} with lifespans $T(\varepsilon)$. Then these solutions satisfy \eqref{Supp_uv}. Moreover, there exists a positive constant $\varepsilon_0=\varepsilon_0(u_0,u_1,v_0,v_1,n,p,q,R)$ such that for any $\varepsilon\in(0,\varepsilon_0]$ the solution $(u,v)$ blows up in finite time. Furthermore, the upper bound estimate for the lifespans
	\begin{align*}
	T(\varepsilon)\leqslant C\varepsilon^{-1/\max\{F_3(n,p,q),F_4(n,p,q)\}}
	\end{align*}
	holds, where $C>0$ is a constant independent of $\varepsilon$ and
	\begin{align}
	F_3(n,p,q)&:= \frac{2+q}{pq-1}-n+1,\label{Defn.F3}\\
	F_4(n,p,q)&:= \frac{1+2p}{pq-1}-n.\label{Defn.F4}
	\end{align}
\end{prop}
Summarizing the above two results, we may immediately obtain the following conclusion for blow-up of energy solutions to Nakao's problem \eqref{Nakao's problem}.
\begin{theorem}\label{Thm.Main.Blow-Up_3}
	Let us consider $p,q>1$ if $n=1,2$, and $1<p,q\leqslant n/(n-2)$ if $n\geqslant 3$ such that
	\begin{align}\label{Assumption.Condition_3}
	\alpha_{\mathrm{N}}:=\max\left\{\frac{q/2+1}{pq-1},\frac{2+p^{-1}}{pq-1},\frac{1/2+p}{pq-1}-\frac{1}{2}\right\}>\frac{n-1}{2}.
	\end{align}
	Let us assume that $(u_0,u_1,v_0,v_1)\in \left(H^1(\mb{R}^n)\times L^2(\mb{R}^n)\right)\times \left(H^1(\mb{R}^n)\times L^2(\mb{R}^n)\right)$ have the same assumption as in Theorem \ref{Thm.Main.Blow-Up}. Let $(u,v)$ be the energy solution to Nakao's problem \eqref{Nakao's problem} according to Definition \ref{DefnEnergySolution} with lifespans $T=T(\varepsilon)$ satisfying \eqref{Supp_uv}. Then, there exists a positive constant $\varepsilon_0=\varepsilon_0(u_0,u_1,v_0,v_1,n,p,q,R)$ such that for any $\varepsilon\in(0,\varepsilon_0]$ the solution $(u,v)$ blows up in finite time $T$. Furthermore, the upper bound estimate for the lifespans
	\begin{align*}
	T(\varepsilon)\leqslant C\varepsilon^{-1/F(n,p,q)}
	\end{align*}
	holds, where $C>0$ is a constant independent of $\varepsilon$ and
	\begin{align*}
	F(n,p,q)&=\max\left\{F_1(n,p,q),F_2(n,p,q),F_3(n,p,q),F_4(n,p,q)\right\}\notag\\
	&=\begin{cases}
	\max\left\{F_3(1,p,q),F_4(1,p,q)\right\}&\mbox{for}\ n=1,\\
	\max\left\{F_1(2,p,q),F_2(2,p,q),F_3(2,p,q),F_4(2,p,q)\right\}&\mbox{for}\ n=2,\\
	\max\left\{F_1(3,p,q),F_4(3,p,q)\right\}&\mbox{for}\ n=3,\\
	F_1(n,p,q)&\mbox{for}\ n\geqslant4.
	\end{cases}
	\end{align*}
\end{theorem}
\begin{remark}
We may observe that the upper bound of lifespan is determined by several components if $n=1,2,3$. Nevertheless, for the higher dimensional cases $n\geqslant 4$, the component $F_1(n,p,q)$ plays a dominant role in the upper bound lifespan estimate we propose.
\end{remark}
\begin{remark}
The conditions $p,q>1$ if $n=1,2$, and $1<p,q\leqslant n/(n-2)$ if $n\geqslant 3$ allow us to guarantee the local (in time) existence of weak solutions to Nakao's problem \eqref{Nakao's problem} with initial data taken from energy space with compact support in a ball with radius $R$. Furthermore, this weak solution $u=u(t,x)$ belongs to classical energy space and has compact support in a ball with radius $R+t$ for $t \in (0,T)$. The proof is given by some standard energy estimates and a contraction mapping argument by using Banach's fixed point theorem.
\end{remark}
\begin{remark}
Let us compare the condition for the exponents of power nonlinearities in our blow-up results, i.e. the conditions \eqref{Assumption.Condition} and \eqref{Assumption.Condition_2}. Under our assumption $p,q>1$ and $p,q\leqslant n/(n-2)$ if $n\geqslant 3$, the condition \eqref{Assumption.Condition_2} partly improves \eqref{Assumption.Condition} only when $n=2$, while for $n\geqslant 3$
\begin{align*}
\left\{(p,q):\alpha_1>\frac{n-1}{2},\ 1<p,q\leqslant\frac{n}{n-2}\right\}\subseteq \left\{(p,q):\alpha_0>\frac{n-1}{2},\ 1<p,q\leqslant\frac{n}{n-2}\right\}.
\end{align*}
\end{remark}
\begin{remark}
 Now, we illustrate the curve $\alpha_{\mathrm{N}}=(n-1)/2$ in different dimensions. Clearly, when $n=1$ the energy solution blows up in finite time for all $1<p,q<\infty$. For this reason, we just describe the case when $n\geqslant 2$ and we will divide the discussion into $n=2$ as well as $n\geqslant 3$, respectively.

\begin{figure}[http]
	\centering
	\begin{tikzpicture}
	\draw[->] (-0.2,0) -- (5.8,0) node[below] {$p$};
	\draw[->] (0,-0.2) -- (0,5.4) node[left] {$q$};
	\node[left] at (0,-0.2) {{$0$}};
	\node[below] at (1,0) {{$1$}};
	\node[left] at (0,1) {{$1$}};
	\node[left, color=red] at (3.72,4.8) {{$\longleftarrow$ $\frac{q/2+1}{pq-1}=\frac{1}{2}$}};
	\node[left, color=blue] at (3.8,4.1) {{$\longleftarrow$ $\frac{2+p^{-1}}{pq-1}=\frac{1}{2}$}};
	\node[left, color=green] at (4.57,3.35) {{$\longleftarrow$ $\frac{1/2+p}{pq-1}-\frac{1}{2}=\frac{1}{2}$}};
	\draw[dashed, color=black]  (0, 1)--(5.6, 1);
	\draw[dashed, color=black]  (1, 0)--(1, 5.4);
	\draw[color=red] plot[smooth, tension=.7] coordinates {(1.2,5.3) (1.8,3.6) (3,2.4)};
	\draw[dashed, color=red] plot[smooth, tension=.7] coordinates {(3,2.4) (4,1.6) (5,1)};
	\draw[color=blue] plot[smooth, tension=.7] coordinates {(3,2.4) (3.5,2.1) (4,2)};
	\draw[dashed, color=blue] plot[smooth, tension=.7] coordinates {(1,5) (1.75,3.4) (3,2.4)};
	\draw[dashed, color=blue] plot[smooth, tension=.7] coordinates {(4,2) (4.8,1.83) (5.4,1.72)};
	\draw[dashed, color=green] plot[smooth, tension=.7] coordinates {(1,4) (2.4,2.4) (4,2)};
	\draw[color=green] plot[smooth, tension=.7] coordinates {(4,2) (5,1.88) (5.4,1.86)};
	\node[left] at (4,-1) {{Case $n=2$}};
	\draw[->] (7.8,0) -- (13.8,0) node[below] {$p$};
	\draw[->] (8,-0.2) -- (8,5.4) node[left] {$q$};
	\node[left] at (8,-0.2) {{$0$}};
	\node[below] at (9,0) {{$1$}};
	\node[left] at (8,1) {{$1$}};
	\node[below] at (12,0) {{$\frac{n}{n-2}$}};
	\node[left] at (8,4) {{$\frac{n}{n-2}$}};
		\node[left, color=red] at (11.72,3.6) {{$\longleftarrow$ $\frac{q/2+1}{pq-1}=\frac{1}{2}$}};
	\node[left, color=blue] at (13.05,2.6) {{$\longleftarrow$ $\frac{2+p^{-1}}{pq-1}=\frac{n-1}{2}$}};
		\node[left, color=green] at (14.3,1.35) {{$\longleftarrow$ $\frac{1/2+p}{pq-1}-\frac{1}{2}=\frac{1}{2}$}};
	\draw[dashed, color=black]  (8, 1)--(13.6, 1);
	\draw[dashed, color=black]  (9, 0)--(9, 5.4);
	\draw[dashed, color=black]  (8, 4)--(13.6, 4);
	\draw[dashed, color=black]  (12, 0)--(12, 5.4);
		\draw[dashed, color=red] plot[smooth, tension=.7] coordinates {(9.2,4) (10,1.6) (11.5,1)};
	\draw[color=blue] plot[smooth, tension=.7] coordinates {(9.5,4) (10.5,2) (12,1.3)};
		\draw[dashed, color=green] plot[smooth, tension=.7] coordinates {(9,3.7) (10,1.8) (12,1.1)};
	\node[left] at (12,-1) {{Case $n\geqslant 3$}};
	\end{tikzpicture}
	\caption{The curve $\alpha_{\mathrm{N}}=(n-1)/2$ in the $p-q$ plane}
	\label{imggg}
\end{figure}
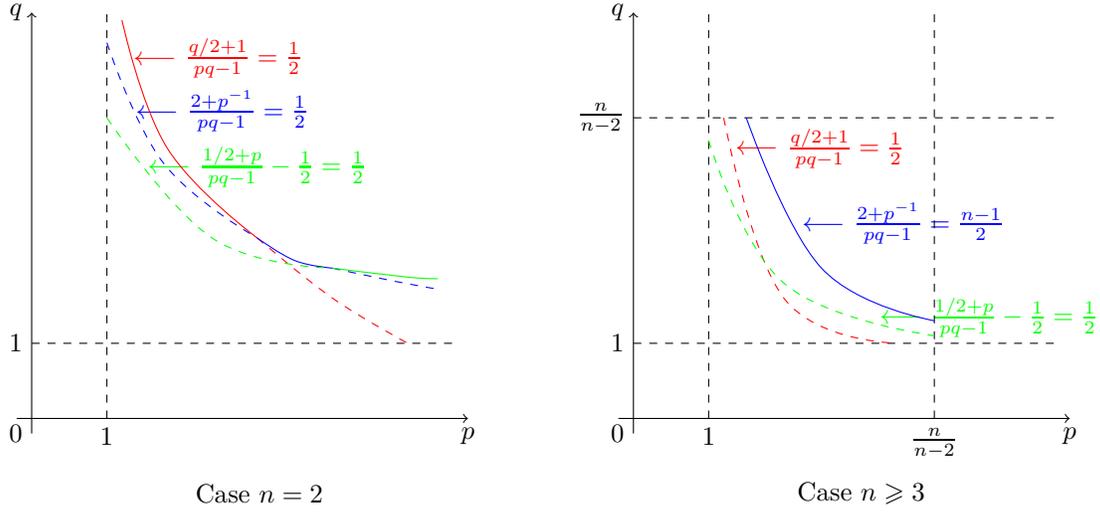
From the graphs, we may observe that in the case when $n=2$, all components in $\alpha_{\mathrm{N}}$ have an influence on the condition \eqref{Assumption.Condition_3}. However, in the case when $n\geqslant3$, under the assumption $1<p,q\leqslant n/(n-2)$, we may derive
\begin{align*}
\alpha_{\mathrm{N}}>\frac{n-1}{2} \, \, \,\, \mbox{iff}\,\,\, \ \  \frac{2+p^{-1}}{pq-1}>\frac{n-1}{2}.
\end{align*}
Consequently, the other components in $\alpha_{\mathrm{N}}$ do not influence our blow-up result.
\end{remark}
\begin{remark}
We may observe from the conditions $1<p,q\leqslant n/(n-2)$ when $n\geqslant 3$ and $\alpha_{\mathrm{N}}>(n-1)/2$ that if $n\geqslant 8$, then energy solutions to Nakao's problem \eqref{Nakao's problem} may blow up for any $1<p,q\leqslant n/(n-2)$.
\end{remark}
\color{black}
\begin{remark}
 We should underline that the critical curve in the $p-q$ plane for Nakao's problem \eqref{Nakao's problem} is still open, especially, the question for global (in time) existence of small data solutions seems to be completely open for $n\geqslant 2$.
\end{remark}
\subsection{Parabolic like versus hyperbolic like} \label{Sec2.1}
Let us turn to some special cases of systems we introduced in Section \ref{Sec1}. If we choose $p=q$ in (\ref{weaklycoupledwave}), then we get
\begin{equation}\label{Eq_p=q_Wave}
\begin{cases}
u_{tt}-\Delta u=|v|^{p},&x\in\mathbb{R}^n,\ t>0,\\
v_{tt}-\Delta v=|u|^{p},&x\in\mathbb{R}^n,\ t>0,\\
(u,u_t,v,v_t)(0,x)=(u_0,u_1,v_0,v_1)(x),&x\in\mathbb{R}^n.
\end{cases}
\end{equation}
Then condition (\ref{Critical curve wave}) reads as follows: \begin{equation*}
\alpha_{\mathrm{W}}:=\frac{p+2+p^{-1}}{p^2-1}=\frac{n-1}{2}.
\end{equation*}
Then, the critical exponent $p_{crit}=p_{crit}(n)$ is the so-called Strauss exponent $p_{\mathrm{Str}}=p_{\mathrm{Str}}(n)$ which is the positive root to the quadratic equation
\begin{align*}
(n-1)p^2-(n+1)p-2=0.
\end{align*}. The above model is hyperbolic like from the point of view of global (in time) existence of small data weak solutions.\\
If we choose $p=q$ in (\ref{weaklycoupleddamped}), then we get
\begin{equation*}
\begin{cases}
u_{tt}-\Delta u+u_t=|v|^{p},&x\in\mathbb{R}^n,\ t>0,\\
v_{tt}-\Delta v+v_t=|u|^{p},&x\in\mathbb{R}^n,\ t>0,\\
(u,u_t,v,v_t)(0,x)=(u_0,u_1,v_0,v_1)(x),&x\in\mathbb{R}^n.
\end{cases}
\end{equation*}
Then condition (\ref{Critical curve damped wave}) reads as follows:
\begin{equation*}
\alpha_{\mathrm{DW}}:=
\frac{1}{p-1}=\frac{n}{2}.
\end{equation*}
Then, the critical exponent $p_{crit}=p_{crit}(n)$ is the so-called Fujita exponent $p_{\mathrm{Fuj}}=p_{\mathrm{Fuj}}(n)=1+\frac{2}{n}$. The above model is parabolic like from the point of view of global (in time) existence of small data weak solutions.

Let us come back to Nakao's problem \eqref{Nakao's problem} carrying $\varepsilon=1$ and $p=q$, namely, to
\begin{equation}\label{Eq_Nakao_p=q}
\begin{cases}
u_{tt}-\Delta u+u_t=|v|^p,&x\in\mb{R}^n, \ t>0,\\
v_{tt}-\Delta v=|u|^p,&x\in\mb{R}^n,\ t>0,\\
(u,u_t,v,v_t)(0,x)=(u_0,u_1,v_0,v_1)(x),&x\in\mb{R}^n.
\end{cases}
\end{equation}
Then, it follows from condition \eqref{Condition Wakasugi} that under the restriction $p\leqslant n/(n-2)$ for $n\geqslant3$ we have blow-up of local (in time) small data weak solutions for
\begin{align}\label{p=q_Wakasugi}
\begin{cases}
1<p<\infty&\mbox{when}\ \ n=1,\\
\displaystyle{1<p\leqslant \max\left\{1+\frac{2}{n},\frac{1}{2(n-1)}\left(1+\sqrt{4n^2-3}\,\right)\right\}}&\mbox{when}\ \ n\geqslant 2.
\end{cases}
\end{align}
Hence, the condition is parabolic like although we have only one semilinear classical damped wave equation in the above model. The influence of the wave equation is missing (at least for large dimensions $n$).

Concerning the condition \eqref{Assumption.Condition_3} stated in Theorem \ref{Thm.Main.Blow-Up_3}, under the restriction $p\leqslant n/(n-2)$ for $n\geqslant3$ we have blow-up of local (in time) small data weak solutions (as long such solutions exist) for
\begin{align}\label{p=q_ChenReis}
\begin{cases}
1<p<\infty&\mbox{when}\ \ n=1,\\
\displaystyle{1<p\leqslant \max\left\{p_0(n),1+\frac{2}{n},\frac{1}{2(n-1)}\left(1+\sqrt{4n^2-3}\,\right)\right\}}&\mbox{when}\ \ n\geqslant 2,
\end{cases}
\end{align}
where $p_0=p_0(n)<p_{\mathrm{Str}}(n)$ for $n \geqslant 2$ is the positive real root of the following cubic equation:
\begin{align*}
(n-1)p^3-(n+3)p-2=0.
\end{align*}

First of all, comparing the condition \eqref{p=q_Wakasugi} with \eqref{p=q_ChenReis}, we observe that an additional component $p_0(n)$ comes in, which means the blow-up range for \eqref{Eq_Nakao_p=q} with the exponent $p$ stated in this paper is larger than those in \cite{Wakasugi2017}. Moreover, the critical exponent for semilinear wave equations \eqref{Eq_p=q_Wave} is the Strauss exponent $p_{\mathrm{Str}}(n)$. 
 In this way we feel that the model is not parabolic like any more,
that is, we feel some influence of the semilinear wave equation.
%
\subsection{Comparison with previous results} \label{Sec2.2}
In this subsection, we will give some remarks and explanations for our result in Theorem \ref{Thm.Main.Blow-Up_3}, especially, the condition
\begin{align*}
\alpha_{\mathrm{N}}=\max\left\{\frac{q/2+1}{pq-1},\frac{2+p^{-1}}{pq-1},\frac{1/2+p}{pq-1}-\frac{1}{2}\right\}>\frac{n-1}{2}
\end{align*}
by comparing with the result from \cite{Wakasugi2017}, in particular, with the condition
\begin{align*}
\alpha_{\mathrm{N,W}}=\max\left\{\frac{q/2+1}{pq-1}+\frac{1}{2},\frac{q+1}{pq-1},\frac{p+1}{pq-1}\right\}\geqslant \frac{n}{2}.
\end{align*}
To guarantee the local (in time) existence of solutions and blow-up in finite time, we assume $p,q>1$ and $p,q\leqslant n/(n-2)$ if $n\geqslant 3$ throughout this subsection.
\begin{itemize}
	\item Concerning the case $n=1$, we prove a blow-up result for energy solutions of Nakao's Cauchy problem when $1<p,q<\infty$. This result corresponds to the blow-up result stated in \cite{Wakasugi2017}. In other words, we may assert that the condition \eqref{Assumption.Condition} is optimal when $n=1$.
	\item Concerning the case $n=2,3$, our result $\alpha_{\mathrm{N}}>(n-1)/2$ partially improves the result in \cite{Wakasugi2017}, that is, the condition $\alpha_{\mathrm{N,W}}\geqslant n/2$. However, some of our results  are worse than those in \cite{Wakasugi2017}. For this reason see the figures below. Particularly, we may divide the condition $\alpha_\mathrm{N}>(n-1)/2$ into three parts in such a way that they are similar to the previous results in the subcritical case as follows:
	\begin{itemize}
		\item Part I: $\frac{q/2+1}{pq-1}+\frac{1}{2}>\frac{n}{2}$, which coincides with the first component in the condition $\alpha_{\mathrm{N,W}}>\frac{n}{2}$;
		\item Part II: $\frac{2+p^{-1}}{pq-1}>\frac{n-1}{2}$, which is similar to the second component in the condition $\alpha_{\mathrm{W}}>\frac{n-1}{2}$;
		\item Part III: $\frac{1/2+p}{pq-1}>\frac{n}{2}$, which is similar to the first component in the condition $\alpha_{\mathrm{DW}}>\frac{n}{2}$.
	\end{itemize}	Moreover, the condition $\alpha_{\mathrm{N,W}}\geqslant n/2$ can be rewritten in the following way: a pair of exponents $(p,q)$ satisfies
	\begin{align*}
	\frac{q/2+1}{pq-1}+\frac{1}{2}\geqslant\frac{n}{2}\ \ \mbox{or} \ \ \alpha_{\mathrm{DW}}\geqslant\frac{n}{2},
	\end{align*}
	whose first part is almost the same as in Part I  (the difference is the limit case). Nevertheless, the main difference between these two conditions are the second and third parts of the above statement.
	\item Concerning the case $n\geqslant 4$, our result completely improves those stated in \cite{Wakasugi2017}. Again, one may see the following figures.
\end{itemize}
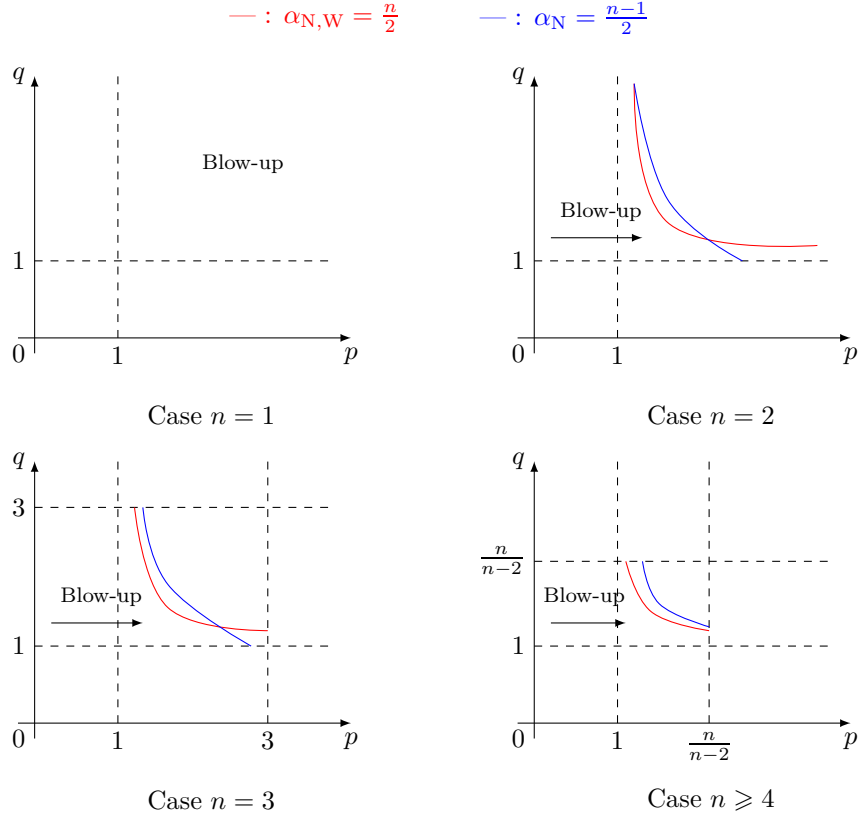
\begin{figure}[http]
	\centering
	\begin{tikzpicture}[>=latex,xscale=1.35,yscale=1.25,scale=0.82]
	\node[left] at (4.5,4.15) {{\color{red} --- : $\alpha_{\mathrm{N,W}}=\frac{n}{2}$}};
	\node[left] at (7.5,4.15) {{\color{blue} --- : $\alpha_{\mathrm{N}}=\frac{n-1}{2}$}};
	\draw[->] (-0.2,0) -- (3.8,0) node[below] {$p$};
	\draw[->] (0,-0.2) -- (0,3.4) node[left] {$q$};
	\node[left] at (0,-0.2) {{$0$}};
	\node[below] at (2.5,2.5) {\footnotesize{Blow-up}};
	\node[below] at (1,0) {{$1$}};
	\node[left] at (0,1) {{$1$}};
	\draw[dashed, color=black]  (0, 1)--(3.6, 1);
	\draw[dashed, color=black]  (1, 0)--(1, 3.4);
	\draw[->] (5.8,0) -- (9.8,0) node[below] {$p$};
	\draw[->] (6,-0.2) -- (6,3.4) node[left] {$q$};
	\node[left] at (6,-0.2) {{$0$}};
	\node[below] at (7,0) {{$1$}};
	\node[left] at (6,1) {{$1$}};
	\draw[dashed, color=black]  (6, 1)--(9.6, 1);
	\draw[dashed, color=black]  (7, 0)--(7, 3.4);
	\draw[color=red] plot[smooth, tension=.7] coordinates {(7.2,3.3) (7.6,1.5) (9.4,1.2)};
	\draw[color=blue]  plot[smooth, tension=.7] coordinates {(7.2,3.3) (7.6,1.8) (8.5,1)};
	\draw[->] (6.2,1.3) -- (7.3,1.3)node[above] at (6.8,1.4) {\footnotesize{Blow-up}};
	\node[left] at (3,-1) {{Case $n=1$}};
	\node[left] at (9,-1) {{Case $n=2$}};
	\draw[->] (-0.2,-5) -- (3.8,-5) node[below] {$p$};
	\draw[->] (0,-5.2) -- (0,-1.6) node[left] {$q$};
	\node[left] at (0,-5.2) {{$0$}};
	\node[below] at (1,-5) {{$1$}};
	\node[left] at (0,-4) {{$1$}};
	\node[below] at (2.8,-5) {{$3$}};
	\node[left] at (0,-2.2) {{$3$}};
	\draw[dashed, color=black]  (0, -4)--(3.6, -4);
	\draw[dashed, color=black]  (1, -5)--(1, -1.6);
	\draw[dashed, color=black]  (0, -2.2)--(3.6, -2.2);
	\draw[dashed, color=black]  (2.8, -5)--(2.8, -1.6);
	\draw[color=red] plot[smooth, tension=.7] coordinates {(1.2,-2.2) (1.6,-3.5) (2.8,-3.8)};
	\draw[color=blue]  plot[smooth, tension=.7] coordinates {(1.3,-2.2) (1.6,-3.2) (2.6,-4)};
	\draw[->] (0.2,-3.7) -- (1.3,-3.7)node[above] at (0.8,-3.6) {\footnotesize{Blow-up}};
	\draw[->] (5.8,-5) -- (9.8,-5) node[below] {$p$};
	\draw[->] (6,-5.2) -- (6,-1.6) node[left] {$q$};
	\node[left] at (6,-5.2) {{$0$}};
	\node[below] at (7,-5) {{$1$}};
	\node[left] at (6,-4) {{$1$}};
	\node[below] at (8.1,-5) {{$\frac{n}{n-2}$}};
	\node[left] at (6,-2.9) {{$\frac{n}{n-2}$}};
	\draw[dashed, color=black]  (6, -4)--(9.6, -4);
	\draw[dashed, color=black]  (7, -5)--(7, -1.6);
	\draw[dashed, color=black]  (6, -2.9)--(9.6, -2.9);
	\draw[dashed, color=black]  (8.1, -5)--(8.1, -1.6);
	\draw[color=red] plot[smooth, tension=.7] coordinates {(7.1,-2.9) (7.4,-3.55) (8.1,-3.8)};
	\draw[color=blue]  plot[smooth, tension=.7] coordinates {(7.3,-2.9) (7.5,-3.45) (8.1,-3.75)};
	\draw[->] (6.2,-3.7) -- (7.1,-3.7)node[above] at (6.6,-3.6) {\footnotesize{Blow-up}};
	\node[left] at (3,-6) {{Case $n=3$}};
	\node[left] at (9,-6) {{Case $n\geqslant4$}};
	\end{tikzpicture}
	\caption{Blow-up range in the $p-q$ plane}
	\label{imgg}
\end{figure}
\medskip

\noindent\textbf{Notation: } We give some notations to be used in this paper. We write $f\lesssim g$ when there exists a positive constant $C$ such that $f\leqslant Cg$. Moreover, $B_R$ denotes the ball around the origin with radius $R$ in $\mathbb{R}^n$.

\section{Proof of Proposition \ref{Thm.Main.Blow-Up}} \label{Sec3}
\setcounter{equation}{0}
\subsection{Iteration frame} \label{Sec3.1}
To begin with the proof, we introduce the following time-dependent functionals:
\begin{align*}
U(t):= \int_{\mathbb{R}^n} u(t,x) \, \mathrm{d} x\ \  \text{and}\ \  V(t):= \int_{\mb{R}^n}v(t,x)\,\mathrm{d}x.
\end{align*}
According to the property of finite propagation speed for weak solutions to wave equations and damped wave equations we know that if the exponents $p,q$ satisfy $1<p,q<\infty$ when $n=1,2,$ and $1<p,q\leqslant n/(n-2)$ when $n\geqslant 3$, and the initial data $(u_0,u_1,v_0,v_1)\in(H^1(\mb{R}^n)\times L^2(\mb{R}^n))^2$ has compact support in $B_R$, then the local (in time) weak solutions belong to the energy space and have compact support in $B_{R+t}$.

Let us choose test functions $\phi$ and $\psi$ in \eqref{Defn.Energy.Solution 01} and \eqref{Defn.Energy.Solution 02} such that $\phi\equiv1$ and $\psi\equiv1$ on $\{(s,x)\in [0,t]\times\mathbb{R}^n: |x|\leqslant R+ s\}$. So, we may immediately derive
\begin{align*}
\int_0^t\int_{\mb{R}^n}u_t(s,x)\,\mathrm{d}x\,\mathrm{d}s+\int_{\mb{R}^n}u_t(t,x)\,\mathrm{d}x-\varepsilon\int_{\mb{R}^n}u_1(x)\,\mathrm{d}x&=\int_0^t\int_{\mb{R}^n}|v(s,x)|^p\,\mathrm{d}x\,\mathrm{d}s,\\
\int_{\mb{R}^n}v_t(t,x)\,\mathrm{d}x-\varepsilon\int_{\mb{R}^n}v_1(x)\,\mathrm{d}x&=\int_0^t\int_{\mb{R}^n}|u(s,x)|^q\,\mathrm{d}x\,\mathrm{d}s.
\end{align*}
The last equations  can be rewritten as
\begin{align}
U'(t)+U(t)&=U'(0)+U(0)+\int_0^t\int_{\mb{R}^n}|v(s,x)|^p\,\mathrm{d}x\,\mathrm{d}s,\label{Pre.Frame.U}\\
V'(t)&=V'(0)+\int_0^t\int_{\mb{R}^n}|u(s,x)|^q\,\mathrm{d}x\,\mathrm{d}s.\label{Pre.Frame.V}
\end{align}
Clearly from \eqref{Pre.Frame.U} and \eqref{Pre.Frame.V}, the functionals fulfill $U(t)\geqslant 0$ and $V(t)\geqslant0$ for any $t\geqslant0$, where we used our assumptions for nonnegative initial data $u_0,u_1,v_0,v_1$. Precisely, \eqref{Pre.Frame.U} gives
\begin{align*}
\frac{\mathrm{d}}{\mathrm{d}t}\left(\mathrm{e}^tU(t)\right)\geqslant \mathrm{e}^t\left(U'(0)+U(0)\right).
\end{align*}
After integration over $[0,t]$, one gets
\begin{align}\label{Lower_Second_U}
U(t)&\geqslant \mathrm{e}^{-t}U(0)+\left(1-\mathrm{e}^{-t}\right)\left(U'(0)+U(0)\right)\notag\\
&=U(0)+\left(1-\mathrm{e}^{-t}\right)U'(0)\geqslant U(0)\geqslant 0.
\end{align}
Concerning the property $V(t)\geqslant0$ for any $t\geqslant0$, we just need to integrate \eqref{Pre.Frame.V} once over $[0,t]$. Thanks to the property of finite propagation speed of $v$ and H\"older's inequality, after multiplying both sides of \eqref{Pre.Frame.U} by $\mathrm{e}^t$ and taking the integration with respect to $t$ over $[0,t]$ we obtain
\begin{align}\label{Iteration Frame U}
U(t)\geqslant\int_0^t\mathrm{e}^{\tau-t}\int_0^{\tau}\int_{\mb{R}^n}|v(s,x)|^p\,\mathrm{d}x\,\mathrm{d}s\,\mathrm{d}\tau\geqslant C_0 \int_0^t\mathrm{e}^{\tau-t}\int_0^{\tau}(R+s)^{-n(p-1)}(V(s))^p\,\mathrm{d}s\,\mathrm{d}\tau
\end{align}
with a positive constant $C_0>0$.  Moreover, taking integration of \eqref{Pre.Frame.V} over $[0,t]$ and considering the finite propagation speed of $u$ show
\begin{align}\label{Iteration Frame V}
V(t)\geqslant \int_0^t\int_0^\tau\int_{\mb{R}^n}|u(s,x)|^q\,\mathrm{d}x\,\mathrm{d}s\,\mathrm{d}\tau\geqslant C_0\int_0^t\int_0^{\tau}(R+s)^{-n(q-1)}(U(s))^q\,\mathrm{d}s\,\mathrm{d}\tau.
\end{align}
All in all, the iteration frames are constructed in \eqref{Iteration Frame U} and \eqref{Iteration Frame V}.
\subsection{Lower bound for the functionals} \label{Sec3.2}
The main approach of the proof is based on an iteration procedure, which requires iteration frames and first lower bound estimates for the functionals $U$ and $V$, respectively. To begin with deriving the estimates of functionals from the below, let us introduce the eigenfunction $\Phi$ of the Laplace operator in $\mb{R}^n$, namely,
\begin{align*}
\Phi(x) & :=  \mathrm{e}^{x}+\mathrm{e}^{-x}  \qquad  \qquad \ \mbox{if} \ n=1, \\  \Phi(x) & := \int_{\mathbb{S}^{n-1}} \mathrm{e}^{x\cdot \omega} \, \mathrm{d} \sigma_\omega \qquad \,\mbox{if} \ n\geqslant 2,
\end{align*}
where $\mathbb{S}^{n-1}$ is the $n-1$ dimensional sphere. Particularly, the function $\Psi$ has been introduced in \cite{YordanovZhang2006}. It satisfies $\Delta\Phi=\Phi$ and has the asymptotic behavior
\begin{align*}
\Phi(x)\sim|x|^{-\frac{n-1}{2}}\,\mathrm{e}^{|x|}\ \ \mbox{as}\ \  |x|\rightarrow\infty.
\end{align*}
By defining the test function with separate variables such that
\begin{align*}
\Psi(t,x):= \mathrm{e}^{-t}\,\Phi(x),
\end{align*}
obviously, the function $\Psi=\Psi(t,x)$ solves the wave equation $\Psi_{tt}-\Delta\Psi=0$.\\
Let us derive lower bound estimates for the functional $U$ in the first place by defining the auxiliary functional
\begin{align*}
V_1(t):=\int_{\mb{R}^n}v(t,x)\,\Psi(t,x)\,\mathrm{d}x.
\end{align*}
According to  \cite{P.Thesis18}, by our nontrivial assumption on initial data satisfying $v_0\not\equiv0$, there exists a constant $C_1=C_1(v_0,v_1)>0$ such that
\begin{align*}
V_1(t)\gtrsim \frac{1}{2}\left(1-\mathrm{e}^{-2t}\right)\int_{\mb{R}^n}(v_0(x)+v_1(x))\,\Phi(x)\,\mathrm{d}x+\mathrm{e}^{-2t}\int_{\mb{R}^n}v_0(x)\,\Phi(x)\,\mathrm{d}x\geqslant C_1\varepsilon
\end{align*}
for any $t\geqslant0$. Indeed, by the asymptotic behavior of the test function $\Psi$, the next inequality holds (see, for example, estimate (2.5) in \cite{LaiTakamura19}):
\begin{align*}
\int_{|x|\leqslant R+t}|\Psi(t,x)|^{\frac{p}{p-1}}\,\mathrm{d}x\leqslant C_2(R+t)^{\frac{(n-1)(2-p')}{2}},
\end{align*}
where $C_2=C_2(n,R)>0$. The application of H\"older's inequality indicates
\begin{align}\label{Eq.2.01}
\int_{\mb{R}^n}|v(t,x)|^p\,\mathrm{d}x\geqslant(V_1(t))^p\bigg(\int_{|x|\leqslant R+t}|\Psi(t,x)|^{\frac{p}{p-1}}\,\mathrm{d}x\bigg)^{-(p-1)}\geqslant C_3\,\varepsilon^p(R+t)^{n-1-\frac{(n-1)p}{2}},
\end{align}
where $C_3:= C_1^pC_2^{1-p}>0$.
By plugging \eqref{Eq.2.01} into \eqref{Iteration Frame U}, we may derive
\begin{align}\label{First Lower bound U}
U(t)&\geqslant C_3\,\varepsilon^p\int_0^t\mathrm{e}^{\tau-t}\int_0^{\tau}(R+s)^{n-1-\frac{(n-1)p}{2}}\,\mathrm{d}s\,\mathrm{d}\tau\notag\\
&\geqslant C_3\,\varepsilon^p(R+t)^{-\frac{(n-1)p}{2}}\int_0^t\mathrm{e}^{\tau-t}\int_0^{\tau}s^{n-1}\,\mathrm{d}s\,\mathrm{d}\tau\notag\\
&\geqslant \frac{C_3\,\varepsilon^p}{n}(R+t)^{-\frac{(n-1)p}{2}}\int_{t/2}^t\mathrm{e}^{\tau-t}\,\tau^n\,\mathrm{d}\tau\notag\\
&= \frac{C_3\,\varepsilon^p}{n2^n}(R+t)^{-\frac{(n-1)p}{2}}\,t^n\left(1-\mathrm{e}^{-t/2}\right)\notag\\
&\geqslant \frac{C_3(1-\mathrm{e}^{-1/2})\varepsilon^p}{n2^n}(R+t)^{-\frac{(n-1)p}{2}}\,t^n
\end{align}
for any $t\geqslant 1$, where $[t/2,t]\subset[0,t]$ has been used.\\
Let us now take the consideration of lower bound estimates for $V$. According to \eqref{Pre.Frame.V} and the nontrivial and nonnegative assumption on $v_1$, the lower bound for the functional $V$ is given by
\begin{align}\label{First Lower bound V}
V(t)\geqslant V(0)+V'(0)t\geqslant C_4\,\varepsilon t
\end{align}
for any $t\geqslant 1$, with a positive constant $C_4=C_4(v_0,v_1)$.

In conclusion, we have obtained first lower bound estimates for the functionals
\begin{align}
U(t)&\geqslant D_1(R+t)^{-\alpha_1}\,t^{\beta_1}\ \mbox{for any}\ t\geqslant1,\label{Initial_U}\\
V(t)&\geqslant Q_1(R+t)^{-a_1}\,t^{b_1}\ \,\mbox{for any}\ t\geqslant1\label{Initial_V},
\end{align}
carrying the multiplicative constants
\begin{align*}
D_1:= \frac{C_3(1-\mathrm{e}^{-1/2})\varepsilon^p}{n2^n}, \ \  Q_1:= C_4\,\varepsilon,
\end{align*}
and the exponents
\begin{align*}
\alpha_1:=\frac{(n-1)p}{2},\ \ a_1:=0,\ \ \beta_1:= n,\ \  b_1:= 1.
\end{align*}
\subsection{Iteration argument} \label{Sec3.3}
In this subsection, we will derive sequences of lower bounds for each functional by using the iteration frames \eqref{Iteration Frame U} and \eqref{Iteration Frame V}. We remark that the iteration argument has been used for weakly coupled systems, for examples, \cite{AKT00,KTW12,PalTak19,PalTak19dt,PalTak19mix}. More precisely, we will show
\begin{align}
U(t)&\geqslant D_j(R+t)^{-\alpha_j}\,(t-L_j)^{\beta_j}\ \mbox{for any}\ t\geqslant L_j,\label{Seq.U(t)}\\
V(t)&\geqslant Q_j(R+t)^{-a_j}\,(t-L_j)^{b_j}\,\ \mbox{for any}\ t\geqslant L_j,\label{Seq.V(t)}
\end{align}
where $\{D_j\}_{j\geqslant 1}$, $\{Q_j\}_{j\geqslant 1}$, $\{\alpha_j\}_{j\geqslant 1}$, $\{a_j\}_{j\geqslant 1}$, $\{\beta_j\}_{j\geqslant 1}$ and $\{b_j\}_{j\geqslant 1}$ are sequences of nonnegative real numbers that will be determined later in the iteration procedure. Moreover, motivated by the recent paper \cite{ChenPalmieri201901}, we construct $\{L_j\}_{j\geqslant 1}$ to be the sequence of the partial products of the convergent infinite product
\begin{align}\label{Seq_Nakao}
\prod\limits_{k=1}^{\infty}\ell_{k}\ \ \mbox{with}\ \ \ell_k:=1+(pq)^{(-k+1)/2}\ \ \mbox{for any}\ \  k\geqslant1,
\end{align}
so that
\begin{align}\label{Seq_Nakao_2}
L_j:=\prod\limits_{k=1}^{j}\ell_{k}\ \ \mbox{for any}\ \  j\geqslant1.
\end{align}
Here, we used the facts that
\begin{align*}
\prod\limits_{k=1}^{\infty}\ell_{k}=\exp\left(\sum\limits_{k=1}^{\infty}\ln\ell_k\right)
\end{align*}
and the ratio test for determining a convergence series
\begin{align*}
\lim\limits_{k\rightarrow\infty}\frac{\ln\ell_{k+1}}{\ln\ell_k}=\lim\limits_{k\rightarrow\infty}\frac{\ln(1+(pq)^{-k/2})}{\ln(1+(pq)^{(-k+1)/2})}=\lim\limits_{k\rightarrow\infty}\frac{(1+(pq)^{(-k+1)/2})(pq)^{-k/2}}{(1+(pq)^{-k})(pq)^{(-k+1)/2}}=(pq)^{-1/2}<1.
\end{align*}
Particularly, the estimates \eqref{Seq.U(t)} and \eqref{Seq.V(t)} have been proved when $j=1$ in the last subsection.

Thus, in order to prove \eqref{Seq.U(t)} and \eqref{Seq.V(t)} by applying an inductive argument, it just remains to show the induction step. On the one hand, let us plug \eqref{Seq.U(t)} into \eqref{Iteration Frame V} and shrink the domain $[0,t]$ into $[L_j,t]$ to get
\begin{align*}
V(t)&\geqslant C_0 D_j^q\int_0^t\int_0^{\tau}(R+s)^{-n(q-1)-q\alpha_j}\,(s-L_j)^{q\beta_j}\,\mathrm{d}s\,\mathrm{d}\tau\\
&\geqslant C_0D_j^q(R+t)^{-n(q-1)-q\alpha_j}\int_{L_j}^t\int_{L_j}^{\tau}(s-L_j)^{q\beta_j}\,\mathrm{d}s\,\mathrm{d}\tau\\
&\geqslant\frac{C_0D_j^q}{(q\beta_j+1)(q\beta_j+2)}(R+t)^{-n(q-1)-q\alpha_j}(t-L_{j+1})^{q\beta_j+2}
\end{align*}
for $t\geqslant L_{j+1}$, where we used $L_{j+1}>L_j$ in the last line of the chain inequality.\\
On the other hand, we combine \eqref{Seq.V(t)} with \eqref{Iteration Frame U} and shrink the domain again. It shows that
\begin{align*}
U(t)&\geqslant C_0Q_j^{p}\int_0^t\mathrm{e}^{\tau-t}\int_0^{\tau}(R+s)^{-n(p-1)-pa_j}\, (s-L_j)^{pb_j}\,\mathrm{d}s\,\mathrm{d}\tau\\
&\geqslant C_0Q_j^{p}(R+t)^{-n(p-1)-pa_j}\int_{L_j}^t\mathrm{e}^{\tau-t}\int_{L_j}^{\tau}(s-L_j)^{pb_j}\,\mathrm{d}s\,\mathrm{d}\tau\\
&\geqslant\frac{C_0Q_j^p}{pb_j+1}(R+t)^{-n(p-1)-pa_j}\int_{t/\ell_{j+1}}^t\mathrm{e}^{\tau-t}\,(\tau-L_j)^{pb_j+1}\,\mathrm{d}\tau\\
&\geqslant\frac{C_0Q_j^p}{pb_j+1}(R+t)^{-n(p-1)-pa_j}(t/\ell_{j+1}-L_j)^{pb_j+1}\int_{t/\ell_{j+1}}^t\mathrm{e}^{\tau-t}\,\mathrm{d}\tau\\
&=\frac{C_0Q_j^p}{(pb_j+1)\ell_{j+1}^{pb_j+1}}(R+t)^{-n(p-1)-pa_j}\,(t-L_{j+1})^{pb_j+1}\left(1-\mathrm{e}^{t(1/\ell_{j+1}-1)}\right)
\end{align*}
for any $t\geqslant L_{j+1}$, which implies $L_j\geqslant t/\ell_{j+1}$. Furthermore, due to the fact that
\begin{align*}
t\geqslant L_{j+1}=\ell_{j+1}L_j\geqslant \ell_{j+1}>1,
\end{align*}
we may estimate
\begin{align*}
1-\mathrm{e}^{t(1/\ell_{j+1}-1)}&\geqslant 1-\mathrm{e}^{-(\ell_{j+1}-1)}\geqslant 1-\left(1-(\ell_{j+1}-1)+\frac{(\ell_{j+1}-1)^2}{2}\right)\\
&=(\ell_{j+1}-1)\left(1-\frac{\ell_{j+1}-1}{2}\right)=(pq)^{-j/2}\left(1-\frac{1}{2(pq)^{j/2}}\right)\\
&=(pq)^{-j}\left((pq)^{j/2}-1/2\right)\geqslant \left((pq)^{1/2}-1/2\right)(pq)^{-j}.
\end{align*}
In other words, we obtain
\begin{align*}
U(t)\geqslant \frac{C_0Q_j^p((pq)^{1/2}-1/2)(pq)^{-j}}{(pb_j+1)\ell_{j+1}^{pb_j+1}}(R+t)^{-n(p-1)-pa_j}\,(t-L_{j+1})^{pb_j+1}
\end{align*}
for any $t\geqslant L_{j+1}$. Thus, \eqref{Seq.U(t)} and \eqref{Seq.V(t)} are valid if the following recursive relations
\begin{align*}
&D_{j+1}:= \frac{C_0Q_j^p((pq)^{1/2}-1/2)(pq)^{-j}}{(pb_j+1)\ell_{j+1}^{pb_j+1}},\ \ \alpha_{j+1}:= n(p-1)+pa_j,\ \ \beta_{j+1}:= pb_j+1,\\
&Q_{j+1}:= \frac{C_0D_j^q}{(q\beta_j+1)(q\beta_j+2)},\ \ \ \ \  \ \ \ \ \ \ \ \ a_{j+1}:= n(q-1)+q\alpha_j,\ \ b_{j+1}:= q\beta_{j}+2,
\end{align*}
are fulfilled.

\subsection{Upper bound estimate for the lifespan}\label{Sub.Sub.Upperbound} \label{Sec3.4}
In the last subsection, we determined the sequence of lower bound estimates for $U$ and $V$. Thus, we may show that the $j$-dependent lower bounds for $U$ and $V$ blows up in finite time when $j\rightarrow\infty$. Simultaneously, the blow-up result and upper bound estimates for the lifespan will be derived.

Let us first determine the explicit representations of $\alpha_j,\beta_j,a_j,b_j$, which contribute to the determination of estimates for the multiplicative constants $D_j$ and $Q_j$.

Concerning the representations of $\alpha_j$ and $a_j$, we just discuss the case when $j$ is an odd integer. For the remaining case that $j$ is an even number, it is unnecessary for the proof of the theorem. By employing the previous definitions for the exponents $\alpha_j$ and $a_j$, one has
\begin{align}
\alpha_j&=n(p-1)+pa_{j-1}=n(pq-1)+pq\alpha_{j-2}=n(pq-1)\sum\limits_{k=0}^{(j-3)/2}(pq)^k+(pq)^{\frac{j-1}{2}}\alpha_1\notag\\
&=(n+\alpha_1)(pq)^{\frac{j-1}{2}}-n=\left(n+\frac{(n-1)p}{2}\right)(pq)^{\frac{j-1}{2}}-n,\label{alphaj}
\end{align}
and by the same approach,
\begin{align}
a_j=(n+a_1)(pq)^{\frac{j-1}{2}}-n=n(pq)^{\frac{j-1}{2}}-n.\label{aj}
\end{align}

Let us now consider the explicit formulas and upper bound estimates for $\beta_j$ and $b_j$ for all $j\geqslant 1$. Combining the definitions of these exponents with an odd number $j$, we claim that
\begin{align}
\beta_j&=pb_{j-1}+1=pq\beta_{j-2}+2p+1=(2p+1)\sum\limits_{k=0}^{(j-3)/2}(pq)^k+(pq)^{\frac{j-1}{2}}\beta_1\notag\\
&=\left(\frac{2p+1}{pq-1}+\beta_1\right)(pq)^{\frac{j-1}{2}}-\frac{2p+1}{pq-1}=\left(\frac{2p+1}{pq-1}+n\right)(pq)^{\frac{j-1}{2}}-\frac{2p+1}{pq-1},\label{betaj}
\end{align}
and similarly,
\begin{align}
b_j&=q\beta_{j-1}+2=pqb_{j-2}+q+2=(q+2)\sum\limits_{k=0}^{(j-3)/2}(pq)^k+(pq)^{\frac{j-1}{2}}b_1\notag\\
&=\left(\frac{q+2}{pq-1}+b_1\right)(pq)^{\frac{j-1}{2}}-\frac{q+2}{pq-1}=\left(\frac{q+2}{pq-1}+1\right)(pq)^{\frac{j-1}{2}}-\frac{q+2}{pq-1}.\label{bj}
\end{align}
In the case when $j$ is an even number (i.e. $j-1$ is an odd number), by the formulas stated in \eqref{betaj} and \eqref{bj}, we may see from the definitions of $\beta_j$ and $b_j$ again that
\begin{align*}
\beta_j&=pb_{j-1}+1=q^{-1}\left(\frac{q+2}{pq-1}+1\right)(pq)^{\frac{j}{2}}-\frac{2p+1}{pq-1},\\
b_j&=q\beta_{j-1}+2=p^{-1}\left(\frac{2p+1}{pq-1}+n\right)(pq)^{\frac{j}{2}}-\frac{q+2}{pq-1}.
\end{align*}
Summarizing the derived representations for odd and even number $j\geqslant 1$, one obtains
\begin{align*}
&\beta_j\leqslant B_0(pq)^{\frac{j-1}{2}}\ \ \mbox{and}\ \  b_j\leqslant\widetilde{B}_0(pq)^{\frac{j-1}{2}}\ \ \mbox{for odd number}\ j,\\
&\beta_j\leqslant B_0(pq)^{\frac{j}{2}}\ \ \,\,\,\,\,\mbox{and}\ \  b_j\leqslant\widetilde{B}_0(pq)^{\frac{j}{2}}\,\ \,\,\,\ \mbox{for even number}\,\,j,
\end{align*}
where $B_0=B_0(p,q,n)$ and $\widetilde{B}_0=\widetilde{B}_0(p,q,n)$ are positive and independent of $j$ constants.

Our next aim is to derive some estimates for $D_j$ and $Q_j$ from the below. It is obviously that
\begin{align*}
pb_{j-1}+1&=\beta_j\leqslant B_0(pq)^{\frac{j}{2}},\\
(q\beta_{j-1}+1)(q\beta_{j-1}+2)&\leqslant (q\beta_{j-1}+2)^2=b_j^2\leqslant\widetilde{B}_0^2(pq)^{j}.
\end{align*}
Moreover, it holds by the application of L'H\^opital's rule
\begin{align*}
\lim\limits_{j\rightarrow\infty}\ell_j^{pb_{j-1}+1}=\lim\limits_{j\rightarrow\infty}\ell_j^{\beta_j}\leqslant\lim\limits_{j\rightarrow\infty}\exp\left(B_0\,(pq)^{j/2}\log\left(1+(pq)^{-j/2}\right)\right)=\exp\left(B_0(pq)^{1/2}\right),
\end{align*}
thus, we may find a suitable constant $M=M(n,p,q)$ such that $\ell_{j}^{-\beta_j}\geqslant M$ for any $j\geqslant 1$. Hence, we can give the following form for lower bounds:
\begin{align*}
D_{j}&=\frac{C_0\,((pq)^{1/2}-1/2)\,(pq)^{-j+1}}{(pb_{j-1}+1)\,\ell_{j}^{pb_{j-1}+1}}Q_{j-1}^p\geqslant \frac{C_0M((pq)^{1/2}-1/2)}{B_0}(pq)^{-\frac{3j}{2}+1}\,Q_{j-1}^p,\\
Q_{j}&=\frac{C_0}{(q\beta_{j-1}+1)(q\beta_{j-1}+2)}D_{j-1}^q\geqslant\frac{C_0}{\widetilde{B}_0^2}(pq)^{-j}\,D_{j-1}^q.
\end{align*}
The derived inequalities immediately lead to
\begin{align}
D_{j}&\geqslant\frac{C_0^{p+1}M((pq)^{1/2}-1/2)}{B_0\widetilde{B}_0^{2p}}(pq)^{-\frac{3j}{2}-(j-1)p+1}\,D_{j-2}^{pq}:= E_0(pq)^{-\frac{3j}{2}-(j-1)p+1}\,D_{j-2}^{pq},\label{Seq.Dj}\\
Q_{j}&\geqslant\frac{C_0^{q+1}M^q((pq)^{1/2}-1/2)^q}{B_0^q\widetilde{B}_0^2}(pq)^{-j-\frac{3(j-1)q}{2}+q}\,Q_{j-2}^{pq}:=\widetilde{E}_0(pq)^{-j-\frac{3(j-1)q}{2}+q}\,Q_{j-2}^{pq}.\label{Seq.Qj}
\end{align}
Considering \eqref{Seq.Dj} with an odd number $j$, we deduce
\begin{align*}
\log D_j&\geqslant pq\log D_{j-2}-\left(\left(\frac{3}{2}+p\right)j-(p+1)\right)\log (pq)+\log E_0\\
&\geqslant (pq)^2\log D_{j-4}-\left(\frac{3}{2}+p\right)(j+(j-2)pq)\log(pq)\\
&\quad+(p+1)(1+pq)\log(pq)+(1+pq)\log E_0\\
&\geqslant (pq)^{\frac{j-1}{2}}\log D_1-\left(\frac{3}{2}+p\right)\log(pq)\sum\limits_{k=1}^{(j-1)/2}\left((j+2-2k)(pq)^{k-1}\right)\\
&\quad+(p+1)\log(pq)\sum\limits_{k=1}^{(j-1)/2}(pq)^{k-1}+\log E_0 \sum\limits_{k=1}^{(j-1)/2}(pq)^{k-1}.
\end{align*}
By an inductive argument, the next formula can be derive
\begin{align*}
\sum\limits_{k=1}^{(j-1)/2}(j+2-2k)(pq)^{k-1}=\frac{1}{pq-1}\left(\frac{2pq}{pq-1}\left(\frac{3}{2}(pq)^{\frac{j-1}{2}}-\frac{1}{2}(pq)^{\frac{j-3}{2}}-1\right)-j\right).
\end{align*}
One has
\begin{align*}
\log D_j&\geqslant (pq)^{\frac{j-1}{2}}\left(\log D_1+\frac{\log(pq)}{2(pq-1)^2}\left(1-7pq-4p^2q\right)+\frac{\log E_0}{pq-1}\right)\\
&\quad+\frac{\log(pq)}{pq-1}\left(\left(\frac{3}{2}+p\right)\left(\frac{2pq}{pq-1}+j\right)-(p+1)\right)-\frac{\log E_0}{pq-1}.
\end{align*}
Thus, for an smallest nonnegative odd number satisfying
\begin{align*}
j\geqslant j_0:= \frac{2(p+1)}{3+2p}+\frac{2\log E_0}{(3+2p)\log(pq)}-\frac{2pq}{pq-1},
\end{align*}
the lower bound can be estimated by
\begin{align}\label{LowBoundDj}
\log D_j&\geqslant (pq)^{\frac{j-1}{2}}\left(\log D_1+\frac{\log(pq)}{2(pq-1)^2}\left(1-7pq-4p^2q\right)+\frac{\log E_0}{pq-1}\right)\notag\\
&=(pq)^{\frac{j-1}{2}}\log \left(D_1(pq)^{-(4p^2q+7pq-1)/(2(pq-1)^2)}E_0^{1/(pq-1)}\right)\notag\\
&=(pq)^{\frac{j-1}{2}}\log \left(E_1\varepsilon^p\right)
\end{align}
for a suitable constant $E_1=E_1(n,p,q)$. By the same way, we may show
\begin{align*}
\log Q_j&\geqslant pq\log Q_{j-2}-\left(\left(1+\frac{3}{2}q\right)j-\frac{5}{2}q\right)\log (pq)+\log \widetilde{E}_0\\
&\geqslant (pq)^2\log Q_{j-2}-\left(1+\frac{3}{2}q\right)\left(j+(j-2)pq\right)\log(pq)\\
&\quad+\frac{5}{2}q(1+pq)\log(pq)+(1+pq)\log\widetilde{E}_0\\
&\geqslant (pq)^{\frac{j-1}{2}}\log Q_1-\left(1+\frac{3}{2}q\right)\log(pq)\sum\limits_{k=1}^{(j-1)/2}\left((j+2-2k)(pq)^{k-1}\right)\\
&\quad+\frac{5}{2}q\log(pq)\sum\limits_{k=1}^{(j-1)/2}(pq)^{k-1}+\log \widetilde{E}_0 \sum\limits_{k=1}^{(j-1)/2}(pq)^{k-1}.
\end{align*}
As a consequence, it yields
\begin{align*}
\log Q_j&\geqslant (pq)^{\frac{j-1}{2}}\left(\log Q_1+\frac{\log(pq)}{(pq-1)^2}\left(-2pq^2-3pq-q+1\right)+\frac{\log\widetilde{E}_0}{pq-1}\right)\\
&\quad+\frac{\log(pq)}{pq-1}\left(\left(1+\frac{3}{2}q\right)\left(\frac{2pq}{pq-1}+j\right)-\frac{5}{2}q\right)-\frac{\log\widetilde{E}_0}{pq-1}.
\end{align*}
If for an odd number $j$ we assume
\begin{align*}
j\geqslant j_1:= \frac{5q}{2+3q}+\frac{2\log\widetilde{E}_0}{(2+3q)\log(pq)}-\frac{2pq}{pq-1},
\end{align*}
then the estimate holds
\begin{align}\label{LowBoundQj}
\log Q_j&\geqslant (pq)^{\frac{j-1}{2}}\left(\log Q_1+\frac{\log(pq)}{(pq-1)^2}\left(-2pq^2-3pq-q+1\right)+\frac{\log\widetilde{E}_0}{pq-1}\right)\notag\\
&\geqslant (pq)^{\frac{j-1}{2}}\log(\widetilde{E}_1\varepsilon)
\end{align}
for a suitable constant $\widetilde{E}_1=\widetilde{E}_1(n,p,q)$.

Let us denote
\begin{align*}
L:= \lim\limits_{j\rightarrow\infty}L_j=\prod\limits_{j=1}^{\infty}\ell_j>1.
\end{align*}
Note that thanks to $\ell_j>1$, it holds $L_j\uparrow L$ as $j\rightarrow\infty$. It leads that \eqref{Seq.U(t)} and \eqref{Seq.V(t)} hold for any odd number $j\geqslant1$ and any $t\geqslant L$.

Let us now consider an odd number $j$ such that $j\geqslant \max\{j_0,j_1\}$. Combining with \eqref{Seq.U(t)}, \eqref{alphaj}, \eqref{betaj} and \eqref{LowBoundDj}, we may observe that
\begin{align*}
U(t)&\geqslant \exp\left((pq)^{\frac{j-1}{2}}\log (E_1\varepsilon^p)\right)(R+t)^{-\alpha_j}(t-L)^{\beta_j}\\
&=\exp\left((pq)^{\frac{j-1}{2}} \left(\log(E_1\varepsilon^p)-\left(\frac{(n-1)p}{2}+n\right)\log(R+t)+\left(\frac{2p+1}{pq-1}+n\right)\log(t-L)\right)\right)\\
&\quad\times(R+t)^{n}(t-L)^{-(2p+1)/(pq-1)}
\end{align*}
for any odd number $j\geqslant\max\{j_0,j_1\}$ and any $t\geqslant L$. Considering $t\geqslant\{R,2L\}$, since $R+t\leqslant 2t$ and $t-L\geqslant t/2$, we have the lower bound estimate for the functional $U$ as follows:
\begin{align}\label{L.B.U(t)}
U(t)\geqslant \exp\left((pq)^{\frac{j-1}{2}}\log\left(E_1\varepsilon^p\,2^{-\frac{(n-1)p}{2}-\frac{2p+1}{pq-1}}\,t^{-\frac{(n-1)p}{2}+\frac{2p+1}{pq-1}}\right)\right)(R+t)^{n}(t-L)^{-(2p+1)/(pq-1)}
\end{align}
for any odd number $j\geqslant\max\{j_0,j_1\}$. The exponent for $t$ in \eqref{L.B.U(t)} can be rewritten by
\begin{align*}
-\frac{(n-1)p}{2}+\frac{2p+1}{pq-1}=p\left(\frac{2+p^{-1}}{pq-1}-\frac{n-1}{2}\right)=pF_1(n,p,q),
\end{align*}
where $F_1(n,p,q)$ is defined in \eqref{Defn.F1}. By our assumption \eqref{Assumption.Condition}, i.e. $F_1(n,p,q)>0$, the exponent for $t$ in the exponential term of \eqref{L.B.U(t)} is positive.

In a similar way as the above, we may derive the next inequality for an odd number $j$ fulfilling $j\geqslant \max\{j_0,j_1\}$:
\begin{align}\label{L.B.V(t)}
V(t)\geqslant\exp\left((pq)^{\frac{j-1}{2}}\log\left(\widetilde{E}_1\varepsilon\,2^{-n-\frac{q+2}{pq-1}-1}t^{-n+\frac{q+2}{pq-1}+1}\right)\right)(t+R)^{n}(t-L)^{-(q+2)/(pq-1)}
\end{align}
for any $j\geqslant\max\{j_0,j_1\}$ and any $t\geqslant L$. Recalling the definition \eqref{Defn.F2}, the exponent for $t$ in \eqref{L.B.V(t)} is
\begin{align*}
-n+\frac{q+2}{pq-1}+1=q\left(\frac{1+2q^{-1}}{pq-1}-\frac{n-1}{q}\right)=qF_2(n,p,q).
\end{align*}
Considering $F_2(n,p,q)>0$ coming from our assumption \eqref{Assumption.Condition}, the exponent for $t$ in the exponential term of \eqref{L.B.V(t)} is positive.

In the case when $F_1(n,p,q)>0$, we set $\varepsilon_0=\varepsilon_0(u_0,u_1,v_0,v_1,n,p,q,R)>0$ such that
\begin{align*}
\left(E_1^{-1}2^{\frac{(n-1)p}{2}+\frac{2p+1}{pq-1}}\right)^{1/(pF_1(n,p,q))}:= E_2\geqslant \varepsilon_0^{1/F_1(n,p,q)}.
\end{align*}
Therefore, for $\varepsilon\in(0,\varepsilon_0]$ and $t>E_2\varepsilon^{-1/F_1(n,p,q)}$ carrying $t\geqslant \max\{R,2L\}$, letting $j\rightarrow\infty$ in \eqref{L.B.U(t)}, we may conclude that the lower bound for $U$ blows up. Similarly, in the remaining case when $F_2(n,p,q)>0$, then we can find a positive constant $\varepsilon_0=\varepsilon_0(u_0,u_1,v_0,v_1,n,p,q,R)>0$ such that
\begin{align*}
\left(\widetilde{E}_1^{-1}2^{n+\frac{q+2}{pq-1}+1}\right)^{1/(qF_2(n,p,q))}:= \widetilde{E}_2\geqslant \varepsilon_0^{1/F_2(n,p,q)}.
\end{align*}
Therefore, for $\varepsilon\in(0,\varepsilon_0]$ and $t>\widetilde{E}_2\varepsilon^{-1/F_2(n,p,q)}$ carrying $t\geqslant \max\{R,2L\}$, letting $j\rightarrow\infty$ in \eqref{L.B.V(t)}, we may conclude that the lower bound for $V$ blows up. In conclusion, these statements proved that the energy solution $(u,v)$ is not globally in time defined and, simultaneously, the lifespan of local (in time) of $(u,v)$ can be estimated by
\begin{align*}
T(\varepsilon)\leqslant C\varepsilon^{-1/\max\{F_1(n,p,q),F_2(n,p,q)\}}.
\end{align*}
The proof of the proposition is complete.

\section{Proof of Proposition \ref{Thm.Main.Blow-Up_2}} \label{Sec4}
\setcounter{equation}{0}
In this section, we will sketch the proof of Proposition \ref{Thm.Main.Blow-Up_2} by using the same techniques as shown in those for Proposition \ref{Thm.Main.Blow-Up}. Nevertheless, we now may provide another lower bound estimates for the functionals $U(t)$ and $V(t)$, which are defined in the previous section.

To begin with, from \eqref{Lower_Second_U} and \eqref{First Lower bound V}, we know there exist positive constants $\widetilde{C}_0$ and $\widetilde{C}_1$ relaying on initial data such that
\begin{align*}
U(t)\geqslant \widetilde{C}_0\,\varepsilon \ \ \mbox{and} \ \ V(t)\geqslant \widetilde{C}_1\,\varepsilon t,
\end{align*}
where the nonnegative and nontrivial assumptions on $u_0$ and $u_1$ were used, respectively. Consequently, by employing H\"older's inequality, supports for solution and the above estimates, we may derive
\begin{align*}
\int_{\mb{R}^n}|u(t,x)|^q\,\mathrm{d}x&\geqslant\widetilde{C}_2(R+t)^{-n(q-1)}(U(t))^q\geqslant \widetilde{C}_0^q\widetilde{C}_2\,\varepsilon^q(R+t)^{-n(q-1)},\\
\int_{\mb{R}^n}|v(t,x)|^p\,\mathrm{d}x&\geqslant\widetilde{C}_3(R+t)^{-n(p-1)}(V(t))^p\geqslant \widetilde{C}_1^p\widetilde{C}_3\,\varepsilon^p(R+t)^{-n(p-1)}\,t^p,
\end{align*}
where $\widetilde{C}_2,\widetilde{C}_3>0$ are suitable constants depending on $n,p,q,R$. They lead that from the estimates \eqref{Iteration Frame U} and \eqref{Iteration Frame V} as follows:
\begin{align*}
U(t)&\geqslant \widetilde{C}_1^p\widetilde{C}_3\,\varepsilon^p\int_0^t\mathrm{e}^{\tau-t}\int_0^{\tau}(R+s)^{-n(p-1)}\,s^p\,\mathrm{d}s\,\mathrm{d}\tau\geqslant\frac{\widetilde{C}_1^p\widetilde{C}_3(1-\mathrm{e}^{-1/2})\varepsilon^p}{(p+1)2^{p+1}}(R+t)^{-n(p-1)}\,t^{p+1},\\
V(t)&\geqslant \widetilde{C}_0^q\widetilde{C}_2\,\varepsilon^q\int_0^t\int_0^{\tau}(R+s)^{-n(q-1)}\,\mathrm{d}s\,\mathrm{d}\tau\geqslant \frac{\widetilde{C}_0^q\widetilde{C}_2\,\varepsilon^q}{2} (R+t)^{-n(q-1)}\,t^2,
\end{align*}
for any $t\geqslant1$. Here, we should underline that the improvement for the lower bound estimates of the functionals in comparison with \eqref{First Lower bound U} and \eqref{First Lower bound V} for large time. The lower bound estimate for the functional $U(t)$ is improved in the case when $n=1$ and $n=2$ with $1<p<2$, moreover, the lower bound estimate for the functional $V(t)$ is improved in the case when $1<q<1+1/n$.

 In other words, the first lower bound estimates \eqref{Initial_U}, \eqref{Initial_V} hold, providing that the multiplicative constants satisfy
\begin{align*}
D_1:=\frac{\widetilde{C}_1^p\widetilde{C}_3(1-\mathrm{e}^{-1/2})\varepsilon^p}{(p+1)2^{p+1}},\ \ Q_1:=\frac{\widetilde{C}_0^q\widetilde{C}_2\,\varepsilon^q}{2}
\end{align*}
and the exponents fulfill
\begin{align*}
\alpha_1:=n(p-1),\ \ a_1:=n(q-1),\ \ \beta_1:=p+1, \ \ b_1:=2.
\end{align*}

Then, following the same approach as those for Theorem \ref{Thm.Main.Blow-Up}, one may derive blow-up of solutions when a pair of exponent $(p,q)$ satisfies
\begin{align*}
-(n+\alpha_1)+\frac{2p+1}{pq-1}+\beta_1=p\left(\frac{2+q}{pq-1}-(n-1)\right)=pF_3(n,p,q)>0,
\end{align*}
or
\begin{align*}
-(n+a_1)+\frac{q+2}{pq-1}+b_1=q\left(\frac{1+2p}{pq-1}-n\right)=q F_4(n,p,q)>0,
\end{align*}
where the functions $F_3(n,p,q)$ and $F_4(n,p,q)$ are defined in \eqref{Defn.F3} and \eqref{Defn.F4}, respectively, moreover, we use our assumption \eqref{Assumption.Condition_2} to guarantee the mentioned functions are positive. Furthermore,  the lifespan of local (in time) of $(u,v)$ can be estimated by
\begin{align*}
T(\varepsilon)\leqslant C\varepsilon^{-1/\max\{F_3(n,p,q),F_4(n,p,q)\}}.
\end{align*}
The proof is complete.


\section*{Acknowledgments}
The PhD study of Wenhui Chen are supported by S\"achsisches Landesgraduiertenstipendium.

\end{document}